\documentclass[english,11pt]{article}
\usepackage{amsfonts}
\usepackage[latin1]{inputenc}
\usepackage{amsmath,amsthm,amssymb}
\usepackage{graphicx}
\usepackage{color}

\def\neweq#1{\begin{equation}\label{#1}}
\def\endeq{\end{equation}}

\newtheorem{theorem}{Theorem}[section]
\newtheorem{lem}{Lemma}[section]

\input{tcilatex}
\begin{document}

\title{\textbf{Quasilinear problems with the competition between convex and
concave nonlinearities and variable potentials}}
\author{Dragos-Patru Covei\break \\
%EndAName
{\small \ Constantin Brancusi University of Tg-Jiu, Calea Eroilor, No 30,
Tg-Jiu, Gorj, Romania}\\
{\small E-mail: }patrucovei@yahoo.com}
\date{}
\maketitle

\begin{abstract}
The purpose of this paper is to prove some existence and non-existence
theorems for the nonlinear elliptic problems of the form \ $-\Delta
_{p}u=\lambda k\left( x\right) u^{q}\pm h\left( x\right) u^{\sigma }$ if $%
x\in \Omega $,\ subject to the Dirichlet conditions $u_{{\footnotesize 1}%
}=u_{2}=0$\ on $\partial \Omega $. In the proofs of our results we use the
sub-super solutions method and variational arguments. Related results as
obtained here have been established in [Z. Guo and Z. Zhang, $W^{1,p}$ 
\textit{versus} $C^{1}$ \textit{local minimizers and multiplicity results
for quasilinear elliptic equations}, Journal of Mathematical Analysis and
Applications, Volume 286, Issue 1, Pages 32-50, 1 October 2003.] for the
case $k\left( x\right) =h\left( x\right) =1$. Our results reveal some
interesting behavior of the solutions due to the interaction between
convex-concave nonlinearities and variable potentials.
\end{abstract}

\baselineskip16pt \renewcommand{\theequation}{\arabic{section}.%
\arabic{equation}} \catcode`@=11 \@addtoreset{equation}{section} \catcode%
`@=12



\section{Introduction and the main results}

In this article we study the existence and non-existence of solutions for
the quasilinear elliptic problems ($P_{\lambda }$)$_{\pm }$ of the following
type 
\begin{equation}
-\Delta _{p}u=\lambda k\left( x\right) u^{q}\pm h\left( x\right) u^{\sigma }%
\text{ if }x\in \Omega \text{, }u>0\text{ in }\Omega \text{, }u=0\text{ of \ 
}\partial \Omega   \tag{($P_{\lambda }$)$_{\pm }$}
\end{equation}%
where $\lambda $ is a positive real parameter, $\Omega \subset \mathbb{R}^{N}
$ ($N\geq 2$) is a bounded domain with smooth boundary, $0<q<p-1<\sigma $,
the variable weight functions $k,$ $h\in L^{\infty }\left( \Omega \right) $
satisfy $ess\inf_{x\in \Omega }k\left( x\right) >0$ and $ess\inf_{x\in
\Omega }h\left( x\right) >0$, and $\Delta _{p}u:=\func{div}\left( \left\vert
\nabla u\right\vert ^{p-2}\nabla u\right) $, $1<p<\infty $ stands for the
p-Laplacian operator.

We call a function $u:\Omega \rightarrow \mathbb{R}$ a solution of problems (%
$P_{\lambda }$)$_{\pm }$ if it belongs to the Sobolev space $%
W_{0}^{1,p}\left( \Omega \right) $ and such that

i)\quad $u\geq 0$ a.e. on $\Omega $ and $u>0$ on a subset of $\Omega $ with
positive measure;

ii)\quad for all $\varphi \in W_{0}^{1,p}\left( \Omega \right) $ the
following identity holds%
\begin{equation*}
\dint\limits_{\Omega }\left\vert \nabla u\right\vert ^{p-2}\nabla u\nabla
\varphi dx=\dint\limits_{\Omega }\left( \lambda k\left( x\right) u^{q}\pm
h\left( x\right) u^{\sigma }\right) \varphi dx.
\end{equation*}%
This kind of problems with convex and concave nonlinearities have been
extensively studied and plays a central role in modern mathematical
sciences, in the theory of heat conduction in electrically conduction
materials, in the study of non-Newtonian fluids (see: Allegretto-Huang \cite%
{AW}, Ambrosetti-Brezis-Cerami \cite{ABC}, Brezis-Nirenberg \cite{BN},
Guo-Zhang \cite{G}, Figueiredo-Gossez-Ubilla \cite{FGU} with their
references). The basic work in our direction is the article \cite{G} where
Guo and Zhang have been considered the Dirichled problem 
\begin{equation*}
-\Delta _{p}u=\lambda u^{q}+u^{\sigma }\text{ if }x\in \Omega \text{, }u>0%
\text{ if }x\in \Omega ,\text{ }u=0\text{ if }x\in \partial \Omega ,
\end{equation*}%
where $\lambda $ is a positive parameter, $\Omega \subset \mathbb{R}^{N}$ ($%
N\geq 3$) is a bounded domain with smooth boundary, $0<q<p-1<\sigma <p^{\ast
}-1$ inequality in which $p^{\ast }$ represents for the Sobolev conjugate
exponent of $p$, namely $p^{\ast }:=Np/\left( N-p\right) $ if $1<p<N$ and $%
p^{\ast }:=\infty $ for $p\geq N$. We mention that in the work \cite{G} the
authors have been extended the results of Brezis and Nirenberg \cite{BN}
obtained in the case $p=2$.

Our main goal is to extend the results obtained in \cite{G} to the more
general problems ($P_{\lambda }$)$_{\pm }$ .

The p-laplacian operator arises naturally in various contexts of physics,
for instance, in non-Newtonian fluid theory, the quantity $p$ is a
characteristic of the medium. The case $1<p<2$ corresponds to pseudoplastics
fluids and $p>2$ arises in the consideration of dilatant fluids.

The main results are as follows:

\begin{theorem}
\label{1} Let\textbf{\ }$p>1$. For all $0<q<p-1<\sigma <p^{\ast }-1$ \textit{%
there exists a positive number }$\lambda ^{\ast }$\textit{\ such that }for $%
\lambda \in \left( 0,\lambda ^{\ast }\right) $ the problem ($P_{\lambda }$)$%
_{+}$ has a minimal solution $u\left( \lambda \right) $ which is increasing
with respect to $\lambda $. \ If \ $\lambda =\lambda ^{\ast }$ the problem ($%
P_{\lambda }$)$_{+}$ has a solution. Moreover, problem ($P_{\lambda }$)$_{+}$
does not have any solution if $\lambda >\lambda ^{\ast }$.
\end{theorem}

\begin{theorem}
\label{2} Suppose $0<q<p-1<\sigma <p^{\ast }-1$. Then \textit{there exists a
positive number }$\lambda ^{\ast }$\textit{\ such that }the problem ($%
P_{\lambda }$)$_{-}$ has at least one solution for $\lambda >\lambda ^{\ast
} $. Moreover, the problem ($P_{\lambda }$)$_{-}$ does not have any solution
for $\lambda <\lambda ^{\ast }$.
\end{theorem}

Before we prove the main theorems, we need some additional results.

\section{Preliminary results}

The next result describes a regularity near the boundary for weak solutions
to (($P_{\lambda }$)$_{\pm }$ ) and is developed by Lieberman in more
general form than one presented here. For the interior regularity we advise
the work of Tolksdorf \cite{tolk2} and DiBenedetto \cite{DIB}.

\begin{lem}
\label{lib}(in \cite{LIB} )\bigskip\ Let $\beta ,\Lambda $, $M_{0}$ be
positive constants with $\beta \leq 1$ and let $\Omega $ be a bounded domain
in $\mathbb{R}^{N}$ with $C^{1,\beta }$ boundary. Suppose $b(x,r)$ satisfies
the condition $\left\vert b(x,r)\right\vert \leq \Lambda $ for all $\left(
x,r\right) $ in $\partial \Omega \times \mathbb{[}-M_{0},M_{0}\mathbb{]}$.
If $u$ is a bounded weak solution of the problem 
\begin{equation*}
\Delta _{p}u+b(x,u)=0\text{ for}\ x\in \Omega \text{, }u=0\text{ on }%
\partial \Omega
\end{equation*}%
with $\left\vert u\right\vert \leq M_{0}$ in $\Omega $, then there is a
positive constant $\alpha :=\alpha \left( \alpha ,\Lambda ,N\right) $ such
that $u$ is in $C^{1,\alpha }\left( \overline{\Omega }\right) $. Moreover $%
\left\vert u\right\vert _{1+\alpha }\leq C\left( \alpha ,\Lambda
,M_{0},N,\Omega \right) $.
\end{lem}

We use in the proof the strong maximum principle of Vazquez.

\begin{lem}
\label{vaz}(see \cite{HOPF}) Let $\Omega $\ be a domain in $\mathbb{R}^{N}$($%
N\geq 1$) and $u\in C^{1}(\Omega )$\ such that $\Delta _{p}u\in
L_{loc}^{2}(\Omega )$, $u\geq 0$\ a.e. in $\Omega $, $u\neq 0$, $\Delta
_{p}u\leq \beta (u)$\ a.e. in \ $\Omega $\ with $\beta :[0,\infty
)\rightarrow \mathbb{R}$\ continuous, non-decreasing, $\beta (0)=0$\ and
either $\beta (s)=0$\ for some $s>0$ or $\beta (s)>0$\ for all $s>0$ \ but%
\begin{equation*}
\int_{0}^{1}(j(S))^{-1/p}dS=\infty \text{ where }j(S)=\int_{0}^{S}\beta
(t)dt.\text{ }
\end{equation*}%
Then if $u$\ does not vanish identically on $\Omega $\ it is positive
everywhere in $\Omega $.\bigskip
\end{lem}

The following lemma has been obtained in Sakaguchi.

\begin{lem}
\label{hopf}(see \cite{SS}) Let $\Omega \subset \mathbb{R}^{N}$ ($N\geq 2$)
is a bounded domain with smooth boundary and let $u\in C^{1}(\overline{%
\Omega })\cap W_{0}^{1,p}\left( \Omega \right) $\ satisfy: 
\begin{eqnarray*}
-\Delta _{p}u &\geq &0\text{ }in\text{ }\Omega \text{ (in the weak sense)},
\\
u &>&0\text{ }in\text{ }\Omega \text{ and }u=0\text{ }on\text{ }\partial
\Omega .
\end{eqnarray*}%
Then $\partial u/\partial n<0$\ on $\partial \Omega $ where $n$ denotes the
unit exterior normal vector to $\partial \Omega $.
\end{lem}

The following comparison principle is proved in \cite{SS} (or consult some
ideas of the proof in \cite[Lemma 3.1.]{tolk}).

\begin{lem}
\label{cp}Let $\Omega \subset \mathbb{R}^{N}$ ($N\geq 2$) is a bounded
domain with smooth boundary and let $u,v\in W_{{}}^{1,p}\left( \Omega
\right) $ satisfy $-\Delta _{p}u\leq -\Delta _{p}v$ for $x\in \Omega $, 
\textit{in the weak sense. }If $u\leq v$ on $\partial \Omega $ then $u\leq v$
in $\Omega $.
\end{lem}

We prove Theorem \ref{1} also by the method of sub- and super-solutions. To
describe this method we introduce the problem%
\begin{equation}
-\Delta _{p}u=\lambda k\left( x\right) u^{q}+h\left( x\right) u^{\sigma }%
\text{ for }x\in \Omega \text{, }u=0\text{ on }\partial \Omega ,  \label{ss}
\end{equation}%
where $\Omega ,\lambda ,k,q,h$ and $\sigma $ are as above. We define $%
\underline{u}\in W_{0}^{1,p}\left( \Omega \right) \cap L^{\infty }\left(
\Omega \right) $ to be a sub-solution of (\ref{ss}) if 
\begin{eqnarray*}
-\Delta _{p}\underline{u} &\leq &\lambda k\left( x\right) \underline{u}%
^{q}+h\left( x\right) \underline{u}^{\sigma }\text{ }x\in \Omega ,\text{ (in
the weak sense)} \\
\underline{u} &=&0,
\end{eqnarray*}%
and $\overline{u}\in W_{0}^{1,p}\left( \Omega \right) \cap L^{\infty }\left(
\Omega \right) $ to be a super-solution of (\ref{ss}) if 
\begin{eqnarray*}
-\Delta _{p}\overline{u} &\geq &\lambda k\left( x\right) \overline{u}%
^{q}+h\left( x\right) \overline{u}^{\sigma }\text{ }x\in \Omega ,\text{ (in
the weak sense)} \\
\overline{u} &=&0.
\end{eqnarray*}%
Then the following result holds:

\begin{lem}
\label{subsuper} (see \cite{XC}) Suppose there exist a sub-solution $%
\underline{u}$ and a super-solution $\overline{u}$ of (\ref{ss}) in the
above sense and that $\underline{u}$ $\leq \overline{u}$. Then there exists
a bounded weak solution $u$ of the problem (\ref{ss}) such that $\underline{u%
}$ $\leq u\leq \overline{u}$.
\end{lem}

We finally recall the following Picone's result for the p-Laplacian
developed by Allegretto and Huang.

\begin{lem}
\label{pic} (see \cite{AW}) Let $v>0$, $u\geq 0$ be differentiable. Denote%
\begin{equation*}
R\left( u,v\right) =\left\vert \nabla u\right\vert ^{p}-\nabla \left( \frac{%
u^{p}}{v^{p-1}}\right) \left\vert \nabla v\right\vert ^{p-2}\nabla v.
\end{equation*}%
Then $R\left( u,v\right) \geq 0$ and $R\left( u,v\right) =0$ a.e. $\Omega $
if and only if $\nabla \left( u/v\right) =0$ a.e. $\Omega $, i.e. $u=kv$ for
some constant $k$ in each component of $\Omega $, where $\Omega $ is bounded
or unbounded, or the whole space $\mathbb{R}^{N}$.
\end{lem}

\section{Proof of the Theorem \protect\ref{1}}

Firstly, we prove that there exists $\lambda _{0}>0$ such that for all $%
\lambda \in (0,\lambda _{0}]$ the problem ($P_{\lambda }$)$_{+}$ has a
solution. The argument relies on constructing a sub- and a super-solution
with the properties from Lemma \ref{subsuper}. In order to find a
sub-solution, consider the problem

\begin{equation}
-\Delta _{p}u=\lambda k\left( x\right) u^{q}\text{ if }x\in \Omega \text{, }%
u>0\text{ in }\Omega \text{, }u=0\text{ on }\partial \Omega .  \label{5d}
\end{equation}%
Then, by \cite{diaz}, problem (\ref{5d}) has a unique positive solution $%
w\in W_{0}^{1,p}\left( \Omega \right) \cap L^{\infty }\left( \Omega \right) $
with $\partial w/\partial n<0$\ on $\partial \Omega $. It is not difficult
to prove that the function $\underline{u}:=\varepsilon ^{1/\left( p-1\right)
}w$ is a sub-solution of problem ($P_{\lambda }$)$_{+}$ provided that $%
\varepsilon >0$ is small enough. For this, it suffices to observe that%
\begin{equation*}
\varepsilon \lambda k\left( x\right) w^{q}\leq \lambda k\left( x\right)
\varepsilon ^{q/\left( p-1\right) }w^{q}+h\left( x\right) \varepsilon
^{\sigma /\left( p-1\right) }w^{\sigma }\text{ in }\Omega
\end{equation*}%
which is true for all $\varepsilon \in \left( 0,1\right) $. Let $v\in
W_{0}^{1,p}\left( \Omega \right) \cap L^{\infty }\left( \Omega \right) $ be
the positive solution of%
\begin{equation*}
\left\{ 
\begin{array}{c}
-\Delta _{p}v=1\text{ in }\Omega \\ 
v=0\text{ on }\partial \Omega .%
\end{array}%
\right.
\end{equation*}%
which exists and is unique from \cite[Lemma 2.1.]{G2}. We prove that if $%
\lambda >0$ is small enough then there is $M>0$ such that $\overline{u}%
=M^{1/\left( p-1\right) }v$ is a super-solution of ($P_{\lambda }$)$_{+}$.
Therefore it suffices to show that%
\begin{equation}
M\geq \lambda k\left( x\right) \left[ M^{1/\left( p-1\right) }v\right]
^{q}+h\left( x\right) \left[ M^{1/\left( p-1\right) }v\right] ^{\sigma }.
\label{6d}
\end{equation}%
In the next, we use some notations%
\begin{equation*}
A=\left\Vert k\right\Vert _{L^{\infty }}\cdot \left\Vert v\right\Vert
_{L^{\infty }}^{q}\text{ and }B=\left\Vert h\right\Vert _{L^{\infty }}\cdot
\left\Vert v\right\Vert _{L^{\infty }}^{\sigma }.
\end{equation*}%
Thus by (\ref{6d}), it is enough to show that there is $M>0$ such that%
\begin{equation*}
M\geq \lambda AM^{q/\left( p-1\right) }+BM^{\sigma /\left( p-1\right) }
\end{equation*}%
that is equivalent to%
\begin{equation}
1\geq \lambda AM^{\left( q-p+1\right) /\left( p-1\right) }+BM^{\left( \sigma
-p+1\right) /\left( p-1\right) }.  \label{7d}
\end{equation}%
Consider the following mapping $\left( 0,\infty \right) \ni t\rightarrow
\lambda At^{\left( q-p+1\right) /\left( p-1\right) }+Bt^{\left( \sigma
-p+1\right) /\left( p-1\right) }$. We also note that this function reaches
its minimum value in $t=C\lambda ^{\left( p-1\right) /\left( \sigma
-q\right) }$, where 
\begin{equation*}
C=\left[ AB^{-1}\left( p-1-q\right) \left( \sigma -p+1\right) ^{-1}\right]
^{\left( p-1\right) /\left( \sigma -q\right) }.
\end{equation*}%
Moreover, the global minimum of this mapping is 
\begin{equation*}
\left[ \left( AC^{\left( q-p+1\right) /\left( p-1\right) }+BC^{\left( \sigma
-p+1\right) /\left( p-1\right) }\right) \right] \lambda ^{\left( \sigma
-p+1\right) /\left( \sigma -p\right) }.
\end{equation*}%
This show that condition (\ref{7d}) is fulfilled for all $\lambda \in \left(
0,\lambda _{0}\right] $ and $M=C\lambda ^{\left( p-1\right) /\left( \sigma
-q\right) },$ where $\lambda _{0}$ satisfies%
\begin{equation*}
\left[ \left( AC^{\left( q-p+1\right) /\left( p-1\right) }+BC^{\left( \sigma
-p+1\right) /\left( p-1\right) }\right) \right] \lambda _{0}^{\left( \sigma
-p+1\right) /\left( \sigma -p\right) }=1.
\end{equation*}%
Taking $\varepsilon >0$ possibly smaller, we also note that the comparison
principle announced in Lemma \ref{cp} implies $\varepsilon ^{1/\left(
p-1\right) }w\leq M^{1/\left( p-1\right) }v$. Thus, by Lemma \ref{subsuper}
the problem ($P_{\lambda }$)$_{+}$ has at least one solution $u_{\lambda }$.
Therefore, this solution is a critical point of the functional 
\begin{equation*}
u\longrightarrow \frac{1}{p}\int_{\Omega }\left\vert \nabla u\right\vert
^{p}dx-\frac{\lambda }{q+1}\int_{\Omega }k\left( x\right) \left\vert
u\right\vert ^{q+1}dx-\frac{1}{\sigma +1}\int_{\Omega }h\left( x\right)
\left\vert u\right\vert ^{\sigma +1}dx
\end{equation*}%
in the closed convex set $\left\{ \left. u\in W_{0}^{1,p}\right\vert
\varepsilon ^{1/\left( p-1\right) }w\leq u\leq M^{1/\left( p-1\right)
}v\right\} $.

By choosing 
\begin{equation*}
\lambda ^{\ast }=\sup \left\{ \left. \lambda >0\right\vert \text{ problem }%
(P_{\lambda })_{+}\text{ has a solution}\right\} ,
\end{equation*}%
we have from the definition of $\lambda ^{\ast }$ that problem $(P_{\lambda
})_{+}$ does not have any solution if $\lambda >\lambda ^{\ast }$. In what
follows we claim that $\lambda ^{\ast }$ is finite. Denote 
\begin{equation*}
m:=\min \left\{ ess\inf_{x\in \Omega }k\left( x\right) \text{, }%
ess\inf_{x\in \Omega }h\left( x\right) \right\} \text{.}
\end{equation*}%
Clearly, $m>0$. Let $\lambda ^{\prime }>0$ be such that 
\begin{equation}
m\left( \lambda ^{\prime }t^{q-p+1}+t^{\sigma -p+1}\right) >\lambda _{1}%
\text{ for all }t\geq 0  \label{1.6d}
\end{equation}%
where $\lambda _{1}$ stands for the first eigenvalue of $\left( -\Delta
_{p}\right) $ in $W_{0}^{1,p}\left( \Omega \right) $. Denote by $\varphi
_{1} $ an eigenfunction of the p-Laplacian operator corresponding to $%
\lambda _{1} $. Then $\varphi _{1}\in C^{1,\alpha }\left( \overline{\Omega }%
\right) $ and $\varphi _{1}>0$ in $\Omega $ as a consequence of the strong
maximum principle of Vazquez (Lemma \ref{vaz}). We apply Picone's result,
Lemma \ref{pic}, to the function $\varphi _{1}$ and $u_{\lambda }$. We drop
the parameter $\lambda $ in the function $u_{\lambda }$ and denote $%
u:=u_{\lambda }$. Observe that $\frac{\varphi _{1}^{p}}{u^{p-1}}$\ belongs
to $W_{0}^{1,p}\left( \Omega \right) $ since $u$ is positive in $\Omega $
and has nonzero outward derivative on the boundary because of the Hopf Lemma %
\ref{hopf}.\ Then for all $\lambda >\lambda ^{\prime }$ we have%
\begin{eqnarray*}
0 &\leq &\int_{\Omega }\left\vert \nabla \varphi _{1}\right\vert
^{p}dx-\int_{\Omega }\nabla \left( \frac{\varphi _{1}^{p}}{u^{p-1}}\right)
\left\vert \nabla u\right\vert ^{p-2}\nabla udx \\
&=&\int_{\Omega }\left\vert \nabla \varphi _{1}\right\vert
^{p}dx-\int_{\Omega }\frac{\varphi _{1}^{p}}{u^{p-1}}\Delta _{p}udx \\
&=&\int_{\Omega }\left\vert \nabla \varphi _{1}\right\vert
^{p}dx-\int_{\Omega }\frac{\varphi _{1}^{p}}{u^{p-1}}\left( \lambda k\left(
x\right) u^{q}+h\left( x\right) u^{\sigma }\right) dx \\
&<&\int_{\Omega }\lambda _{1}\varphi _{1}^{p}dx-\int_{\Omega }m\left[
\lambda k\left( x\right) u^{q-p+1}+h\left( x\right) u^{\sigma -p+1}\right]
\varphi _{1}^{p}dx \\
&<&\int_{\Omega }\lambda _{1}\varphi _{1}^{p}dx-\int_{\Omega }m\left[
\lambda ^{\prime }u^{q-p+1}+u^{\sigma -p+1}\right] \varphi _{1}^{p}dx \\
&=&\int_{\Omega }\left[ \lambda _{1}-m\left( \lambda ^{\prime
}u^{q-p+1}+u^{\sigma -p+1}\right) \right] \varphi _{1}^{p}dx<0.
\end{eqnarray*}%
Thus we get a desired contradiction. As a conclusion we obtain the following
result $\lambda ^{\ast }\leq \lambda ^{\prime }<\infty $ which proves our
claim. Let as now prove that $u_{\lambda }$ is a minimal solution of the
problem ($P_{\lambda }$)$_{+}$. By the definition of $\lambda ^{\ast }$
there exists $\overline{\lambda }<\lambda $ such that $\overline{\lambda }%
<\lambda ^{\ast }$ and ($P_{\overline{\lambda }}$)$_{+}$ has a positive
solution $u_{\overline{\lambda }}$. The rest of the argument is based on the
standard monotone iteration. Consider the sequence $\left( u_{n}\right)
_{n\geq 0}$ defined by $u_{0}=w$ (where $w$ is the unique solution of (\ref%
{5d})) and $u_{n}$ the solution of the problem 
\begin{eqnarray*}
-\Delta _{p}u_{n} &=&\lambda k\left( x\right) u_{n-1}^{q}+h\left( x\right)
u_{n-1}^{\sigma },\text{ if }x\in \Omega \\
u_{n}\left( x\right) &>&0,\text{ if }x\in \Omega \\
\text{ }u_{n}\left( x\right) &=&0,\text{ if }x\in \partial \Omega
\end{eqnarray*}%
which exists and is unique from the results in \cite{LL} (see also arguments
in \cite{G}). By using the comparison principle, it is not hard to show that 
\begin{equation}
u_{0}=w\leq u_{1}\leq ...\leq u_{n}\leq u_{n+1}\leq u_{\overline{\lambda }}%
\text{ in }\Omega .  \label{in}
\end{equation}%
In fact, it follows again by the above cited comparison principle applied to
the problem%
\begin{eqnarray*}
-\Delta _{p}u_{0} &=&\lambda k\left( x\right) u_{0}^{q}\leq \lambda k\left(
x\right) u_{0}^{q}+h\left( x\right) u_{0}^{\sigma }=-\Delta _{p}u_{1}\text{
in }\Omega , \\
u_{0} &=&u_{1}=0\text{ on }\partial \Omega
\end{eqnarray*}%
that $u_{0}\leq u_{1\text{ }}$in $\Omega $. Similarly, one can show by using
the same Lemma \ref{cp} that $u_{1}\leq u_{2}$ in $\Omega $. In particular,
for all $x\in \Omega $ the sequence $\left( u_{n}\right) _{n\geq 0}$ is a
nondecreasing sequence which is bounded and therefore $u_{n}\leq U$ for any
positive solution $U$ of ($P_{\lambda }$)$_{+}$. \ Using the relation (\ref%
{in}), the decay property of $u_{\overline{\lambda }}$ and a standard
diagonalization procedure we get a subsequence converging to a solution $%
u_{\lambda }$ of ($P_{\lambda }$)$_{+}$, satisfying $u_{\lambda }\leq u_{%
\overline{\lambda }}$ and $u_{\lambda }\leq U$ for any arbitrary solution $U$
of problem ($P_{\lambda }$)$_{+}$. The conclusion then follow. At this stage
it is easy to deduce that the mapping $u_{\lambda }$ is increasing with
respect to $\lambda $. We consider $u_{\lambda _{1}}$, $u_{\lambda _{2}}$
with $0<\lambda _{1}<\lambda _{2}<\lambda ^{\ast }$. Since%
\begin{equation*}
-\Delta _{p}u_{\lambda _{2}}=\lambda _{2}k\left( x\right) u_{\lambda
_{2}}^{q}+h\left( x\right) u_{\lambda _{2}}^{\sigma }>\lambda _{1}k\left(
x\right) u_{\lambda _{2}}^{q}+h\left( x\right) u_{\lambda _{2}}^{\sigma }
\end{equation*}%
then $u_{\lambda _{2}}$ is a super-solution of problem ($P_{\lambda _{1}}$)$%
_{+}$. The argument used above may be used to construct a sequence $\left(
u_{n}\right) _{n\geq 0}$ such that $0<u_{n-1}<u_{n}<u_{\lambda _{2}}$
converging to a solution $U$ of ($P_{\lambda _{1}}$)$_{+}$ with $%
U<u_{\lambda _{2}}$ and therefore $u_{\lambda _{1}}\leq U<u_{\lambda _{2}}$
by the minimality of $u_{\lambda _{1}}$. This proves our claim.

It remain to show that problem ($P_{\lambda }$)$_{+}$ has a solution if $%
\lambda =\lambda ^{\ast }$. For this purpose it is enough to prove that 
\begin{equation}
u_{\lambda }\text{ is bounded in }W_{0}^{1,p}\left( \Omega \right) \text{ as 
}\lambda \rightarrow \lambda ^{\ast }.  \label{bound}
\end{equation}%
Thus, by passing to a suitable subsequence if necessary, we may assume 
\begin{equation*}
u_{\lambda }\rightarrow u^{\ast }\text{ in }W_{0}^{1,p}\left( \Omega \right) 
\text{ as }\lambda \rightarrow \lambda ^{\ast },
\end{equation*}%
which implies that $u^{\ast }$ is a weak solution of ($P_{\lambda }$)$_{+}$
provided that $\lambda =\lambda ^{\ast }$. Moreover since the mapping $%
\lambda \rightarrow u_{\lambda }$ is increasing, it follows that $u^{\ast
}\geq 0$ a.e. on $\Omega $ and $u^{\ast }>0$ on a subset of $\Omega $ with
positive measure. As we mentioned, it is often advantageous to work with $u$
instead of $u_{\lambda }$. A key ingredient of the proof is that all
solutions $u$ have negative energy. More precisely, if $E:W_{0}^{1,p}\left(
\Omega \right) \rightarrow \mathbb{R}$ is defined by%
\begin{equation*}
E\left( u\right) =\frac{1}{p}\int_{\Omega }\left\vert \nabla u\right\vert
^{p}dx-\frac{\lambda }{q+1}\int_{\Omega }k\left( x\right) \left\vert
u\right\vert ^{q+1}dx-\frac{1}{\sigma +1}\int_{\Omega }h\left( x\right)
\left\vert u\right\vert ^{\sigma +1}dx
\end{equation*}%
then%
\begin{equation}
E\left( u\right) <0\text{ for all }\lambda \in \left( 0,\lambda ^{\ast
}\right) \text{.}  \label{8d}
\end{equation}%
We do it in the following steps:

\textit{Step 1)}\quad the solution $u$ satisfies%
\begin{equation}
\int_{\Omega }\left\{ \left\vert \nabla u\right\vert ^{p}-\left[ \lambda
q/\left( p-1\right) \right] k\left( x\right) u^{q+1}+\left[ \sigma /\left(
p-1\right) \right] h\left( x\right) u^{\sigma +1}\right\} dx\geq 0.
\label{9d}
\end{equation}%
This follows by the same arguments from \cite[Lemma 3.7.]{G}.

\textit{Step 2)}\quad Since $u$ is a solution of ($P_{\lambda }$)$_{+}$ we
have 
\begin{equation}
\int_{\Omega }\left\vert \nabla u\right\vert ^{p}dx=\int_{\Omega }\lambda
k\left( x\right) u^{q+1}dx+\int_{\Omega }h\left( x\right) u^{\sigma +1}dx.
\label{10d}
\end{equation}%
Plugging relation (\ref{9d}) into (\ref{10d}) we have 
\begin{equation}
\lambda \left( p-1-q\right) \int_{\Omega }k\left( x\right) u^{q+1}dx\geq
\left( \sigma +1-p\right) \int_{\Omega }h\left( x\right) u^{\sigma +1}dx
\label{11d}
\end{equation}%
In particular, it follows from these two latest relations that%
\begin{eqnarray*}
E\left( u\right) &=&\lambda \left( \frac{1}{p}-\frac{1}{q+1}\right)
\int_{\Omega }k\left( x\right) u^{q+1}dx+\left( \frac{1}{p}-\frac{1}{\sigma
+1}\right) \int_{\Omega }h\left( x\right) u^{\sigma +1}dx \\
&=&-\lambda \frac{p-1-q}{p\left( q+1\right) }\int_{\Omega }k\left( x\right)
u^{q+1}dx+\frac{\sigma +1-p}{p\left( \sigma +1\right) }\int_{\Omega }h\left(
x\right) u^{\sigma +1}dx \\
&\leq &-\lambda \frac{p-1-q}{p\left( q+1\right) }\int_{\Omega }k\left(
x\right) u^{q+1}dx+\lambda \frac{p-1-q}{p\left( \sigma +1\right) }%
\int_{\Omega }h\left( x\right) u^{\sigma +1}dx\leq 0.
\end{eqnarray*}%
Thus, by combining (\ref{8d}) and (\ref{9d}), sobolev embedings, and using
the fact that $k$, $h\in L^{\infty }\left( \Omega \right) $ it follows 
\begin{equation*}
\sup \left\{ \left. \left\Vert u_{\lambda }\right\Vert _{W_{0}^{1,p}\left(
\Omega \right) }\right\vert \lambda <\lambda ^{\ast }\right\} <\infty
\end{equation*}%
and so (\ref{bound}) is finished. This complete the proof of Theorem \ref{1}.

\section{Proof of the Theorem \protect\ref{2}}

\bigskip The study of existence of solutions to problem ($P_{\lambda }$)$%
_{-} $ is done by looking for critical points of the functional $F_{\lambda
}:W_{0}^{1,p}\left( \Omega \right) \rightarrow \mathbb{R}$ defined by%
\begin{equation*}
F_{\lambda }\left( u\right) =\frac{1}{p}\int_{\Omega }\left\vert \nabla
u\right\vert ^{p}dx-\frac{\lambda }{q+1}\int_{\Omega }k\left( x\right)
\left\vert u\right\vert ^{q+1}dx+\frac{1}{\sigma +1}\int_{\Omega }h\left(
x\right) \left\vert u\right\vert ^{\sigma +1}dx.
\end{equation*}%
In the next we adopt the following notations 
\begin{equation*}
\left\Vert u\right\Vert :=\left( \int_{\Omega }\left\vert \nabla
u\right\vert ^{p}dx\right) ^{1/p}\text{, }\left\Vert u\right\Vert
_{q+1}:=\left( \int_{\Omega }\left\vert u\right\vert ^{q+1}dx\right)
^{1/(q+1)}\text{, }\left\Vert u\right\Vert _{\sigma +1}:=\left( \int_{\Omega
}\left\vert u\right\vert ^{\sigma +1}dx\right) ^{1/(\sigma +1)}.
\end{equation*}%
We prove that $F_{\lambda }$ is coercive. In order to verify this claim, we
first observe that%
\begin{equation*}
F_{\lambda }\left( u\right) \geq \frac{1}{p}\left\Vert u\right\Vert
^{p}-C_{1}\left\Vert u\right\Vert _{q+1}^{q+1}+C_{2}\left\Vert u\right\Vert
_{\sigma +1}^{\sigma +1},
\end{equation*}%
where 
\begin{equation*}
C_{1}=\frac{\lambda }{q+1}\left\Vert k\right\Vert _{L^{\infty }}\text{ and }%
C_{2}=\frac{1}{\sigma +1}ess\inf_{x\in \Omega }h\left( x\right)
\end{equation*}%
are positive constants. Since $q<\sigma $, a simple calculation shows that
the mapping 
\begin{equation*}
\left( 0,\infty \right) \ni t\rightarrow C_{1}t^{q+1}-C_{2}t^{\sigma +1}
\end{equation*}%
attains its global minimum $m<0$ at 
\begin{equation*}
t=\left[ \frac{C_{2}\left( q+1\right) }{C_{1}\left( \sigma +1\right) }\right]
^{1/\left( \sigma -q\right) }.
\end{equation*}%
So we conclude that%
\begin{equation*}
F_{\lambda }\left( u\right) \geq \frac{1}{p}\left\Vert u\right\Vert ^{p}+m%
\text{,}
\end{equation*}%
and hence $F_{\lambda }\left( u\right) \rightarrow \infty $ as $\left\Vert
u\right\Vert \rightarrow \infty $ whish finished the proof that $F_{\lambda
} $ is coercive. Let $\left( u_{n}\right) _{n\geq 0}$ be a minimizing
sequence of $F_{\lambda }$ in $W_{0}^{1,p}\left( \Omega \right) $. The
coercivity of $F_{\lambda }$ implies the boundedness of $u_{n}$ in $%
W_{0}^{1,p}\left( \Omega \right) $. Then, up to a subsequence if necessary,
we may assume that there exists \ $u$\ in $W_{0}^{1,p}\left( \Omega \right) $%
\ non-negative such that $u_{n}\overset{n\rightarrow \infty }{\rightarrow }u$
weakly in $W_{0}^{1,p}\left( \Omega \right) $. We remark that the function \ 
$u$\ can be\ non-negative due to $F_{\lambda }\left( u\right) =F_{\lambda
}\left( \left\vert u\right\vert \right) $. Standard arguments based on the
lower semi-continuity of the energy functional show that $u$ is a global
minimizer of $F_{\lambda }$ and therefore is a solution in the sense of
distributions of ($P_{\lambda }$)$_{-}$.

In what follows we claim that the weak limit $u$ is a non-negative weak
solution of problem ($P_{\lambda }$)$_{-}$ if $\lambda >0$ is large enough.
We first observe that $F_{\lambda }\left( 0\right) =0$. So, in order to
prove that the non-negative solution is non-trivial, it suffices to prove
that there exists $\Lambda >0$ such that%
\begin{equation*}
\inf_{u\in W_{0}^{1,p}\left( \Omega \right) }F_{\lambda }\left( u\right) <0%
\text{ for all }\lambda >\Lambda .
\end{equation*}%
For this purpose we consider the variational problem with constraints,%
\begin{equation}
\Lambda =\inf \left\{ \left. \frac{1}{p}\int_{\Omega }\left\vert \nabla
v\right\vert ^{p}dx+\frac{1}{\sigma +1}\int_{\Omega }h\left( x\right)
\left\vert v\right\vert ^{\sigma +1}dx\right\vert \text{ }v\in
W_{0}^{1,p}\left( \Omega \right) \text{ and }\frac{1}{q+1}\int_{\Omega
}k\left( x\right) \left\vert v\right\vert ^{q+1}dx=1\right\} .  \label{12d}
\end{equation}%
Let $\left( v_{n}\right) _{n\geq 0}$ be an arbitrary minimizing sequence for
this problem. Then $v_{n}$ is bounded, hence we can assume that it weakly
converges to some $v\in W_{0}^{1,p}\left( \Omega \right) $ with%
\begin{equation*}
\frac{1}{q+1}\int_{\Omega }k\left( x\right) \left\vert v\right\vert
^{q+1}dx=1\text{ and }\Lambda =\frac{1}{p}\int_{\Omega }\left\vert \nabla
v\right\vert ^{p}dx+\frac{1}{\sigma +1}\int_{\Omega }h\left( x\right)
\left\vert v\right\vert ^{\sigma +1}dx.
\end{equation*}%
Thus 
\begin{equation*}
F_{\lambda }\left( v\right) =\Lambda -\lambda \text{ for all }\lambda
>\Lambda \text{.}
\end{equation*}%
Set%
\begin{equation*}
\lambda ^{\ast }:=\inf \left\{ \left. \lambda >0\right\vert \text{ problem }%
(P_{\lambda })_{-}\text{ admits a nontrivial weak solution}\right\} \geq 0%
\text{.}
\end{equation*}%
The above remarks show that $\Lambda \geq \lambda ^{\ast }$ and that problem 
$(P_{\lambda })_{-}$ has a solution for all $\lambda \geq \Lambda $. We now
argue that problem $(P_{\lambda })_{-}$ has a solution for all $\lambda
>\lambda ^{\ast }$. Fixed $\lambda >\lambda ^{\ast }$, by the definition of $%
\lambda ^{\ast }$, we can take $\mu \in \left( \lambda ^{\ast }\text{, }%
\lambda \right) $ such that $F_{\mu \text{ }}$ has a nontrivial critical
point $u_{\mu }\in W_{0}^{1,p}\left( \Omega \right) $. Since $\mu <\lambda $%
, it follows that $u_{\mu }$ is a sub-solution of problem $(P_{\lambda
})_{-} $. We now want to construct a super-solution that dominates $u_{\mu }$%
. For this purpose we consider the constrained minimization problem%
\begin{equation}
\inf \left\{ F_{\lambda }\left( v\right) ,\text{ }v\in W_{0}^{1,p}\left(
\Omega \right) \text{ and }v\geq u_{\mu }\right\} .  \label{13d}
\end{equation}%
From the previous arguments, used to treat (\ref{12d}) follows that problem (%
\ref{13d}) has a solution $u_{\lambda }>u_{\mu }$. Moreover, $u_{\lambda }$
is a solution of problem $(P_{\lambda })_{-}$ for all $\lambda >\lambda
^{\ast }$. \ With the arguments developed in \cite{G} we deduce that problem 
$(P_{\lambda })_{-}$ has a solution if $\lambda =\lambda ^{\ast }$. The same
monotonicity arguments as above show that $(P_{\lambda })_{-}$\ does not
have any solution if $\lambda <\lambda ^{\ast }$. Fix $\lambda >\lambda
^{\ast }$. It remains to argue that the non-negative weak solution $u$ is,
in fact, positive. Indeed, using Moser iteration, we obtain that $u\in
L^{\infty }\left( \Omega \right) $. Once $u\in L^{\infty }\left( \Omega
\right) $ it follows by Lemma \ref{lib} that $u$ is a $C^{1,\alpha }\left( 
\overline{\Omega }\right) $ solution of problem $(P_{\lambda })_{-}$
provided for some $\alpha $. Invoking the nonlinear strong maximum principle
of Vazquez (Lemma \ref{vaz}), since $u$ is a non-negative smooth weak
solution of the differential inequality 
\begin{equation*}
-\Delta _{p}u+h\left( x\right) u^{\sigma }\geq 0\text{ in }\Omega ,\text{ }
\end{equation*}%
we deduce that $u$ is positive everywhere in $\Omega $. The proof of Theorem %
\ref{2} is completed.

The extension of the above results to all space $\mathbb{R}^{N}$ or to the
nonlinearities depending on the gradient $\nabla u$ requires some further
nontrivial modifications and will be considered in a future work. We
anticipate that the methods and concepts here can be extended to systems or
when in discussion are more general linear/non-linear operators as well.

\end{document}